\documentstyle{article}

\begin{document}

\hfill{q-alg/9610011}

\hfill{October, 1996}

\vspace{10mm}

\centerline{\bf  QUANTUM ORTHOGONAL CAYLEY-KLEIN GROUPS }
\centerline{ \bf IN CARTESIAN BASIS}
\vspace*{0.37truein}
\centerline{\footnotesize N.A. GROMOV, I.V. KOSTYAKOV and V.V. KURATOV}
\vspace*{0.015truein}
\centerline{\footnotesize\it Department of Mathematics, Komi Science Centre,}
\baselineskip=10pt
\centerline{\footnotesize\it Ural Division, Russian Academy of Sciences,}
\baselineskip=10pt
\centerline{\footnotesize\it Syktyvkar, 167000, Russia }
\baselineskip=10pt
\centerline{\footnotesize\it E-mail : parma@omkomi.intec.ru}
\vspace*{0.21truein}

 \abstract {The similarity transformations of quantum orthogonal
groups are developed and FRT theory is reformulated to the
Cartesian basis. The quantum orthogonal Cayley-Klein groups are
introduced as the algebra functions over an associative algebra
with the nilpotent generators. The  quantum orthogonal
Cayley-Klein algebras are obtained as the dual objects to
the corresponding quantum  groups.}

\textheight=7.8truein
\textwidth=5truein
\setcounter{footnote}{0}
\renewcommand{\thefootnote}{\alph{footnote}}

\section{Introduction}
\noindent
     The {\it simple} quantum groups and algebras are described
by FRT  theory{\cite{1}}.  Although  the   quantum   group   as   a
mathematical structure is not a group at all but it is a Hopf
algebra, nevertheless there are many parallel features for both
group and quantum group.

     In group theory there is a remarkable set of groups, namely
the motion groups of $ 3^n $  $ n $-dimensional spaces of constant
curvature or
the orthogonal Cayley-Klein (CK) groups{\cite{2}}. The classical
simple orthogonal groups $ SO(n), $ the semisimple pseudoorthogonal
groups $ SO(p,q) $ and nonsemisimple groups such as Euclidean
$ E(n), $ Poincare $ P(n), $ Galileian $ G(n) $ are in this set.
The nonsemisimple groups may be obtained from the simple (or
semisimple) one's by the well known procedure of
contraction{\cite{Wig--53}}. On  the  other  hand the set of CK groups
may be described in uniform manner with the help of pure
algebraical methods{\cite{2}}.

     In the present paper we apply this approach for description of
quantum (non\-semi\-simple) orthogonal CK groups and corresponding
quantum algebras. As it turned out slightly modified FRT theory
is suitable for this purpose: it is enough to replace the
complex number field {\bf C} by the dual algebra
$ {\bf D}_{n}(\iota;{\bf C}) $ with the nilpotent generators.

     The contents of this work is as follows.
     The dual algebra $ {\bf D}_{n}(\iota;{\bf C}) $ is briefly
  described in Sec.2. The
 linear transformations of generators
  of quantum  orthogonal  group  are  regarded  in  Sec.3  and  the
description of $ SO_{q}(N;{\bf C}) $ is given. In Sec.4 the quantum
Cayley-Klein groups $ SO_{v}(N;j) $ are introduced as Hopf algebra
of the   noncommutative  functions  with  dual  variables.  Quantum
Cayley-Klein algebra $ so_{v}(N;j) $ are cosidered in Sec.5 as
the dual object to the $ SO_{v}(N;j). $ The developed theory is
illustrated by the N=3 example of  quantum  group  and  algebra  in
Sec.6.

\section{ Dual algebra $ {\bf D}_{n}(\iota;{\bf C}) $ }
\noindent
Dual algebra $ {\bf D}_{n}(\iota;{\bf C}) $ is defined as an associative
  algebra with unit and {\it nilpotent} generators
   $ \iota_{1},\ldots,\iota_{n},\; \iota_{k}^{2}=0,\; k=1,\ldots,n $
 with {\it commutative} multiplication
     $ \iota_{k}\iota_{m}=\iota_{m}\iota_{k},\; k\neq{m}. $
  The general element of $ {\bf D}_{n}(\iota;{\bf C}) $ has the form
\begin{equation}
a=a_0 + \sum_{p=1}^{2^{n}-1} \sum_{k_1<...<k_p}
a_{k_1...k_p}\iota_{k_1}...\iota_{k_p}, \quad
a_0,a_{k_1...k_p}\in{\bf C}.
\end{equation}
For $ n=1 $ we have
$ {\bf D}_{1}(\iota_{1};{\bf C})\ni{a=a_{0}+a_{1}\iota_{1}}, $
i.e. dual (or Study) numbers, when
$ a_{0},a_{1}\in{\bf R}. $
For $ n=2 $ the general element of
$ {\bf D}_{2}(\iota_{1},\iota_{2};{\bf C}) $
is written as follows:
$ a=a_{0}+a_{1}\iota_{1}+a_{2}\iota_{2}+a_{12}\iota_{1}\iota_{2}. $

Divisions $ a/{\iota_k}, a \in {\bf R}, {\bf C},$ and
$ {\iota_m}/{\iota_k}, k \neq m $ are not defined, but division of
a dual unit by itself is equal to real unit
$ {\iota_k}/{\iota_k}=1.$ A function of a dual
argument is defined by its Taylor expansion
$ f(x_0+{\iota_k}x_1)=f(x_0)+{\iota_k}x_1{f^{'}(x_0)}.$

The well known Grassmann algebra
   $ \Gamma_{n}( \xi ) $
 is the algebra with {\it nilpotent} generators
     $ \xi_k^2=0, k=1,\ldots,n $
 and {\it anticommutative} multiplication
     $ \xi_k\xi_m= -\xi_m\xi_k, k \neq m. $
  It is easy to verify that the product of two generators of
 Grassmann algebra has the same algebraic properties as the
 generator of dual algebra
     $ \iota_k=\xi_k\xi_{n+k}, \; k=1,\ldots,n. $
  This means that the dual algebra is the subalgebra of even
  part of Grassmann algebra
     $ {\bf D}_{n}(\iota;{\bf C}) \subset   \Gamma_{2n}(\xi ). $

\section{ Cartesian basis for $ SO_{q}(N;{\bf C}) $ }
\noindent

     According with FRT theory\cite{1} of quantum groups let  us
regard an algebra $ {\bf C} \langle t_{ik} \rangle $
of noncommutative polynomials of
$ N^2 $ variables $ t_{ik}, i,k=1,\ldots,N $ over complex number
field ${\bf C}. $ For well known\cite{1} lower triangular matrix
$ R_{q}\in{M_{N^2}({\bf C})} $ the generators
$ T=(t_{ik})_{i,k=1}^{N}\in{M_{N}({\bf C} \langle t_{ik} \rangle)} $
have the following commutation relations
\begin{equation}
RT_1T_2=T_2T_1R,
\label{2}
\end{equation}
where $ T_1=T \otimes I, \  T_2=I \otimes T \in M_{N^2}({\bf C}
\langle t_{ij} \rangle). $
There are additional relations of $ q$-orthogonality
\begin{equation}
TCT^t=T^tCT=C,
\label{3}
\end{equation}
where $ C=C_0q^{\rho},  \rho=diag(\rho_1, \ldots, \rho_N), $
\begin{equation}
(\rho_1, \ldots, \rho_N)=
\left \{ \begin{array}{ccc}
     (n-\frac{1}{2}, n-\frac{3}{2}, \ldots , \frac{1}{2},0,-\frac{1}{2},
     \ldots , -n+\frac{1}{2}), & {\rm for} & N=2n+1 \\
     (n-1, n-2, \ldots, 1,0,0,-1, \ldots, -n+1), & {\rm for} & N=2n.
     \end{array} \right.
\label{4}
\end{equation}
and the only nonzero elements of matrix
$ C_{0}\in{M_{N}({\bf C})} $
are real units on the second diagonal
     $ (C_{0})_{ik}=\delta_{i'k}, i'=N+1-i. $
The quantum orthogonal group $ SO_{q}(N;{\bf C}) $
is defined as the quotient
     \begin{equation}
SO_{q}(N;{\bf C})={\bf C} \langle t_{ik} \rangle \big/ (\ref{2}),(\ref{3}).
    \label{5}
    \end{equation}
>From the algebraic point of view $ SO_{q}(N;{\bf C}) $ is
a Hopf algebra with the following coproduct $ \Delta, $
counit $ \epsilon $ and antipode $ S: $
\begin{equation}
\Delta T=T \dot {\otimes} T,\quad \epsilon(T)=I,\quad S(T)=CT^tC^{-1}.
\label{6}
\end{equation}

     For $ q=1 $ Eq.(\ref{3}) is the orthogonality condition
written with the help of matrix $ C_{0} $ and the set of
matrices $ T $ is the orthogonal group  $ SO(N;{\bf C}) $
in so-called "symplectic" basis. We shall call $ SO_{q}(N;{\bf C}) $
described by Eqs.(\ref{2})--(\ref{5}) as the orthogonal
quantum group in symplectic basis.

     One of the solutions of matrix equation
\begin{equation}
DC_0D^t=I,
\label{7}
\end{equation}
     have the form
$$
D=\frac{1}{\sqrt{2}}
\left ( \begin{array}{cc}
      I &      {\tilde C_0} \\
      i{\tilde C_0} &   -iI
      \end{array} \right ),    \       N=2n,
$$
\begin{equation}
D=\frac{1}{\sqrt{2}}
\left ( \begin{array}{ccc}
      I & 0 &   {\tilde C_0} \\
      0 & \sqrt{2} &  0 \\
      i{\tilde C_0} & 0  &  -iI
      \end{array} \right ),    \       N=2n+1,
\label{8}
\end{equation}
     where $ \tilde{C}_{0}\in{M_{n}({\bf C})} $ is the matrix with
 the real units on the second diagonal. For the nondeformed
 case $ q=1 $ the equation
     \begin{equation}
   U=DTD^{-1}
    \label{9}
    \end{equation}
define the similarity transformation of the orthogonal matrix
from symplectic to Cartesian basis. So we shall call the linear
transformation (\ref{9}) of generators of quantum group
$ SO_{q}(N;{\bf C}) $ as transformation to the Cartesian basis.

     It is easy to reformulate Eqs.(\ref{2})--(\ref{5}) to
the new basis. Commutation relations of Cartesian (or rotation)
generators
  $ U=(u_{ik})_{i,k=1}^{N}\in{M_{N}({\bf C} \langle u_{ik} \rangle )} $
are given by
\begin{equation}
 \tilde{R_{q}}U_{1}U_{2} = U_{2}U_{1} \tilde{R}_{q},
\label{10}
\end{equation}
where
\begin{equation}
\tilde{R}_{q} = (D \otimes D)R_{q}(D \otimes D)^{-1}
\label{11}
\end{equation}
is now  the  nontriangular  matrix.  The  additional  relations  of
$ q $-orthogonality are as follows
\begin{equation}
UC'U^t=C',\quad
U^t(C')^{-1}U=(C')^{-1},
\label{12}
\end{equation}
  where
\begin{equation}
C'=DCD^t=
=\left ( \begin{array}{ccc}
\cosh{z\tilde \rho} & 0 &i \sinh{z\tilde \rho}\tilde ó_0 \\
0     &      1    &    0 \\
-i{\tilde C_0} \sinh{z\tilde \rho} & 0 & {\tilde C_0} \cosh{z\tilde \rho}
{\tilde ó_0}
\end{array} \right ),
\label{13}
\end{equation}
   $ \tilde{\rho}=diag(\rho_{1},\ldots,\rho_{n})\in{M_{n}({\bf C})}, $
   for $ N=2n+1 $ and without middle column and row for $ N=2n. $

     The orthogonal quantum group $ SO_{q}(N;{\bf C}) $ in Cartesian basis
 is described by the quotient
     \begin{equation}
SO_{q}(N;{\bf C})={\bf C} \langle u_{ik} \rangle \big/ (\ref{10}),(\ref{12})
    \label{14}
    \end{equation}
and is Hopf algebra with the following coproduct $ \Delta, $
counit $ \epsilon $ and antipode $ S: $
\begin{equation}
\Delta U=U \dot \otimes U,\quad \epsilon(U)=I,\quad S(U)=C'U^t(C')^{-1}.
\label{15}
\end{equation}

\section{
     Quantum CK groups $ SO_{v}(N;j) $ as matrix groups
             over dual algebras ${\bf D}_{N-1}(\iota). $ }
\noindent
 Orthogonal Cayley-Klein (CK) groups $ SO(N;j) $ (or motion groups
 of spaces of constant curvature)\cite{2} are realized in the
 Cartesian basis as the matrix groups over ${\bf D}_{N-1}(\iota) $
 with the help of the {\it special} matrices
     \begin{displaymath}
   (A(j))_{kp}=\tilde{J}_{kp}a_{kp}, \quad a_{kp}\in{C},
    \end{displaymath}
     \begin{equation}
  \tilde{J}_{kp}=J_{kp},\; k<p,   \quad    \tilde{J}_{kp}=J_{pk},\;
k\geq{p}, \quad J_{\mu\nu}=\prod_{r=\mu}^{\nu -1}j_{r}, \;
     \mu < \nu, \; j_{r}=1,\iota_{r},i.
    \label{16}
    \end{equation}
 These matrices   are   satisfied   the  following $ j$-orthogonality
relations:
     \begin{equation}
    A(j)A^{t}(j)=A^{t}(j)A(j)=I.
    \label{17}
    \end{equation}

     In the symplectic basis the orthogonal CK groups are described
 by the matrices
     \begin{equation}
 B(j)=D^{-1}A(j)D
    \label{18}
    \end{equation}
 with the additional relations of $ j$-orthogonality
     \begin{equation}
   B(j)C_{0}B^{t}(j)=B^{t}(j)C_{0}B(j)=C_{0}.
    \label{19}
    \end{equation}
The matrix $ D $ in Eqs.(\ref{18}) is given by Eqs.(\ref{8}).

   For example matrices from CK group $ SO(3;j) $ in Cartesian
     basis are
     \begin{equation}
        A(j) =
     \left( \begin{array}{ccc}
     a_{11}      & j_{1}a_{12}  & j_{1}j_{2}a_{13} \\
j_{1}a_{21}      & a_{22}       & j_{2}a_{23}      \\
j_{1}j_{2}a_{31} & j_{2}a_{32}  & a_{33} \\
     \end{array} \right)
    \label{19'}
    \end{equation}
 and in symplectic basis are
     \begin{equation}
        B(j) =
     \left( \begin{array}{ccc}
 b_{11}+ij_{1}j_{2}\tilde{b}_{11} & j_{1}b_{12}-ij_{2}\tilde{b}_{12}
& b_{13}-ij_{1}j_{2}\tilde{b}_{13} \\
j_{1}b_{21}+ij_{2}\tilde{b}_{21}  & b_{22}
& j_{1}b_{21}-ij_{2}\tilde{b}_{21}       \\
 b_{13}+ij_{1}j_{2}\tilde{b}_{13} & j_{1}b_{12}+ij_{2}\tilde{b}_{12}
& b_{11}-ij_{1}j_{2}\tilde{b}_{11} \\
     \end{array} \right) .
    \label{19''}
    \end{equation}

     We shall regard the quantum deformations of the contracted CK
 groups, i.e. $ j_{k}=1,\iota_{k}. $ We shall follow to FRT
theory\cite{1}, but   instead  of  the  algebra  $  {\bf  C}\langle
t_{ik}\rangle $ our starting
 point is  $ {\bf D}\langle t_{ik} \rangle  $  ---  the  algebra   of
   noncommutative
 polynomials of $ N^2 $ variables $ t_{ik},\; i,k=1,\ldots,N $  over
 dual algebra $ {\bf D}_{N-1}(\iota). $ In addition we shall transform
 the deformation parameter $ q=\exp{z} $ as follows:
     \begin{equation}
     z=Jv, \quad J\equiv J_{1N}=\prod_{k=1}^{N-1}j_{k},
    \label{20}
    \end{equation}
 where $ v $ is the new deformation parameter.

     {\it In symplectic basis} the quantum CK group $ SO_{v}(N;j) $
 is produced by the generating matrix
$ T(j)\in{M_{N}({\bf D} \langle t_{ik} \rangle ) } $
 equal to B(j) (\ref{18}) for $ q=1. $ The non\-com\-mu\-ta\-ti\-ve entries of
 $ T(j) $ obey the commutation relations
\begin{equation}
 R_v(j)T_1(j)T_2(j)=T_2(j)T_1(j)R_v(j).
\label{21}
\end{equation}
 and the additional relations of $ (v,j)$-orthogonality
     \begin{equation}
   T(j)C(j)T^{t}(j)=T^{t}(j)C(j)T(j)=C(j),
    \label{22}
    \end{equation}
 where lower triangular R--matrix $ R_{v}(j) $ and $ C(j) $ are
 obtained from $ R_q $ and $ C, $ respectively, by substitution
 $ Jv $ instead of $ z: $
     \begin{equation}
     R_{v}(j)=R_{q}(z \rightarrow Jv), \quad
     C(j)=C(z \rightarrow Jv).
    \label{23}
    \end{equation}
 Then the quotient
     \begin{equation}
SO_{v}(N;j)= {\bf D} \langle t_{ik} \rangle \big/ (\ref{21}),(\ref{22})
    \label{24}
    \end{equation}
is Hopf algebra with the  coproduct $ \Delta, $
counit $ \epsilon $ and antipode $ S: $
\begin{equation}
\Delta T(j)=T(j) \dot {\otimes}T(j),\quad \epsilon (T(j))=I, \quad
 S(T(j))=C(j)T^t(j)C^{-1}(j).
\label{25}
\end{equation}

     {\it In Cartesian basis} the quantum CK group  $ SO_{v}(N;j) $
 is generated    by   $   U(j)=(\tilde{J}_{ik}u_{ik})\in{M_{N}({\bf
D}\langle u_{ik} \rangle) } $
 with the commutation relations
\begin{equation}
\tilde R_v(j)U_1(j)U_2(j)=U_2(j)U_1(j)\tilde R_v(j)
\label{26}
\end{equation}
 and additional relations of $ (v,j)$-orthogonality
\begin{equation}
U(j)C'(j)U^t(j)=C'(j), \quad
U^t(j)(C'(j))^{-1}U(j)=(C'(j))^{-1},
\label{27}
\end{equation}
 where
     \begin{equation}
 \tilde{R}_{v}(j)=\tilde{R}_{q}(z \rightarrow Jv), \quad
     C'(j)=C'(z \rightarrow Jv).
    \label{28}
    \end{equation}
 The quotient
     \begin{equation}
SO_{v}(N;j)={\bf D}\langle u_{ik} \rangle \big/ (\ref{26}),(\ref{27})
    \label{28'}
    \end{equation}
is Hopf algebra with the following  coproduct $ \Delta, $
counit $ \epsilon $ and antipode $ S: $
\begin{equation}
\Delta U(j)=U(j) \dot {\otimes}U(j),\quad \epsilon (U(j))=I, \quad
 S(U(j))=C^{'}(j)U^t(j)(C^{'}(j))^{-1}.
\label{29}
\end{equation}

\section{ Quantum CK algebras $ so_{v}(N;j) $ as a dual to
     $ SO_{v}(N;j). $ }

\noindent

 According to FRT theory\cite{1} the dual space
   $ Hom(SO_{v}(N;j),{\bf C}) $ is an algebra with the multiplication
 induced by coproduct $ \Delta $ in $ SO_{v}(N;j) $
     \begin{equation}
     l_{1}l_{2}(a)=(l_{1}\otimes l_{2})(\Delta (a)),
    \label{30}
    \end{equation}
  $ l_{1},l_{2}\in{ Hom(SO_{v}(N;j),{\bf C})}, \quad a\in{SO_{v}(N;j)}. $
 Let us formally introduce $ N \times N $ upper $ (+) $ and lower
 $ (-) $ triangular matrices $ L^{(\pm)}(j) $ as follows: it is
 necessary to put $ j_{k}^{-1} $ in the nondiagonal matrix elements
  of $ L^{(\pm)}(j), $ if there is the parameter $ j_k $ in the
  corresponding matrix element of $ T(j). $ For example, if
  $ (T(j))_{12}=j_{1}t_{12}+j_{2}\tilde{t}_{12}, $ then
  $ (L^{(+)}(j))_{12}=j_{1}^{-1}l_{12}+j_{2}^{-1}\tilde{l}_{12}. $
  Formally the matrices   $   L^{(\pm)}(j)   $   are   not   defined  for
  $ j_{k}=\iota_{k}, $ since $ \iota_{k}^{-1} $ do not exist, but
  if we give an action of the matrix functionals $ L^{(\pm)}(j) $
  on the elements of $ SO_{v}(N;j) $ by the duality relation
     \begin{equation}
    \langle L^{(\pm)}(j),T(j) \rangle = R^{(\pm)}(j),
    \label{31}
    \end{equation}
where
     \begin{equation}
     R^{(+)}(j)=PR_{v}(j)P, \quad  R^{(-)}(j)= R_{v}^{-1}(j),\quad
     Pu \otimes w = w \otimes u,
    \label{32}
    \end{equation}
then we shall have well defined expressions even for dual values
of the parameters $ j_{k}. $

     The elements of $ L^{(\pm)}(j) $ satisfy the commutation
 relations
     \begin{eqnarray}
     R^{(+)}(j)L_{1}^{(\sigma)}(j)L_{2}^{(\sigma)}(j) & = &
 L_{2}^{(\sigma)}(j)L_{1}^{(\sigma)}(j)R^{(+)}(j), \nonumber  \\
     R^{(+)}(j)L_{1}^{(+)}(j)L_{2}^{(-)}(j) & = &
L_{2}^{(-)}(j)L_{1}^{(+)}(j)R^{(+)}(j),  \quad
     \sigma = \pm
    \label{33}
    \end{eqnarray}
and additional relations
     \begin{eqnarray}
  L^{(\pm)}(j)C^{t}(j)L^{(\pm)}(j) & = & C^{t}(j),  \nonumber \\
L^{(\pm)}(j)(C^{t}(j))^{-1}L^{(\pm)}(j)  & = & (C^{t}(j))^{-1}, \nonumber \\
l_{kk}^{(+)}l_{kk}^{(-)}=l_{kk}^{(-)}l_{kk}^{(+)}=1, & &
     l_{11}^{(+)}\ldots l_{NN}^{(+)}=1, \; k=1,\ldots ,N.
    \label{34}
    \end{eqnarray}
 An algebra  $  so_{v}(N;j)=\{I,L^{(\pm)}(j)\} $ is called quantum CK
algebra and
is Hopf algebra with the following  coproduct $ \Delta, $
counit $ \epsilon $ and antipode $ S: $
  $$
 \Delta L^{(\pm)}(j)=L^{(\pm)}(j) \dot {\otimes} L^{(\pm)}(j), \quad
\epsilon (L^{(\pm)}(j))=I,
  $$
\begin{equation}
S(L^{(\pm)}(j))  =  C^{t}(j)(L^{(\pm)}(j))^{t}(C^{t}(j))^{-1}.
    \label{35}
    \end{equation}

 It is possible to show that algebra $ so_{v}(N;j) $ is isomorphic
 with the quantum deformation of the universal enveloping algebra
 of the CK algebra $ so(N;j), $ which may be obtained from the
 orthogonal algebra $ so(N) $ by contractions\cite{2}.
So there are at least two ways for constraction of quantum CK
algebras.
\vspace{5mm}

\section{ Example: $ SO_{v}(3;j) $ and $ so_{v}(3;j) $
in symplectic basis }
\noindent

The generating matrix for the simplest quantum orthogonal group
  $ SO_{v}(3;j),\; j=(j_1,j_2) $ is in the form
     \begin{equation}
        T(j) =
     \left( \begin{array}{ccc}
 t_{11}+ij_{1}j_{2}\tilde{t}_{11} & j_{1}t_{12}-ij_{2}\tilde{t}_{12}
& t_{13}-ij_{1}j_{2}\tilde{t}_{13} \\
j_{1}t_{21}+ij_{2}\tilde{t}_{21}  & t_{22}
& j_{1}t_{21}-ij_{2}\tilde{t}_{21}       \\
 t_{13}+ij_{1}j_{2}\tilde{t}_{13} & j_{1}t_{12}+ij_{2}\tilde{t}_{12}
& t_{11}-ij_{1}j_{2}\tilde{t}_{11} \\
     \end{array} \right) .
    \label{36}
    \end{equation}
The R-matrix is obtained from the standart one by Eq.(\ref{23})
and is as follows
$$
        R_{v}(j)\equiv R_{q}(z \rightarrow Jv) =
$$
     \begin{equation}
     \left( \begin{array}{ccccccccc}
e^{Jv}  & 0 & 0 & 0 & 0 & 0 & 0 & 0 & 0  \\
     0  & 1 & 0 & 0 & 0 & 0 & 0 & 0 & 0  \\
     0  & 0 & e^{-Jv} & 0 & 0 & 0 & 0 & 0 & 0  \\
     0  & 2\sinh Jv  & 0 & 1 & 0 & 0 & 0 & 0 & 0  \\
          0  & 0 & -2e^{-Jv/2}\sinh Jv & 0 & 1 & 0 & 0 & 0 & 0  \\
          0  & 0 & 0 & 0 & 0 & 1 & 0 & 0 & 0  \\
 0 & 0 & 2(1-e^{-Jv})\sinh Jv & 0 & -2e^{-Jv/2}\sinh Jv & 0 & e^{-Jv} &0& 0 \\
          0  & 0 & 0 & 0 & 0 & 2\sinh Jv & 0 & 1 & 0  \\
          0  & 0 & 0 & 0 & 0 & 0 & 0 & 0 & e^{Jv}  \\
     \end{array} \right) .
    \label{39}
    \end{equation}
 Qver the dual algebras $ {\bf D}_2(j_1,j_2),\; j_1=\iota_1, j_2=1, $ or $
j_1=1, j_2=\iota_2, $ or $  j_1=\iota_1, j_2=\iota_2 $ this R-matrix
  may be written in the form
     \begin{equation}
   R_{v}(j) = I + Jv \tilde R,
    \label{40}
    \end{equation}
where
     \begin{equation}
     \tilde R =
     \left( \begin{array}{ccccccccc}
1  & 0 & 0 & 0 & 0 & 0 & 0 & 0 & 0  \\
     0  & 0 & 0 & 0 & 0 & 0 & 0 & 0 & 0  \\
     0  & 0 & -1 & 0 & 0 & 0 & 0 & 0 & 0  \\
     0  & 2  & 0 & 0 & 0 & 0 & 0 & 0 & 0  \\
          0  & 0 & -2 & 0 & 0 & 0 & 0 & 0 & 0  \\
          0  & 0 & 0 & 0 & 0 & 0 & 0 & 0 & 0  \\
          0  & 0 & 0 & 0 & -2 & 0 & -1  & 0 & 0  \\
          0  & 0 & 0 & 0 & 0 & 2 & 0 & 0 & 0  \\
          0  & 0 & 0 & 0 & 0 & 0 & 0 & 0 & 1  \\
     \end{array} \right)  .
    \label{41}
    \end{equation}
The commutation  relations  and  additional  relations  of  $ (v,j)
$-orthogonality may be obtained from Eqs.(\ref{21}),(\ref{22}) by
straitforward calculations, so we shall concentrate our attention
on the constraction of quantum algebra $ so_{v}(3;j). $

     The matrix functionals $ L^{\pm}(j) $ have the form
     \begin{equation}
        L^{(+)}(j) =
     \left( \begin{array}{ccc}
 l_{11} & j_{1}^{-1}l_{12}-ij_{2}^{-1}\tilde{l}_{12}
& l_{13}-ij_{1}^{-1}j_{2}^{-1}\tilde{l}_{13} \\
0  & 1 & j_{1}^{-1}l_{21}-ij_{2}^{-1}\tilde{l}_{21}       \\
 0 & 0 & l_{11}^{-1} \\
     \end{array} \right) ,
    \label{37}
    \end{equation}
     \begin{equation}
        L^{(-)}(j) =
     \left( \begin{array}{ccc}
l_{11}^{-1}  & 0 & 0 \\
j_{1}^{-1}l_{21}+ij_{2}^{-1}\tilde{l}_{21}  & 1 & 0       \\
l_{13}+ij_{1}^{-1}j_{2}^{-1}\tilde{l}_{13}  & j_{1}^{-1}l_{12}+ij_{2}^{-1}
\tilde{l}_{12}
&  l_{11}  \\
     \end{array} \right) .
    \label{38}
    \end{equation}
Their actions on the generators (\ref{36}) of quantum group
  $ SO_{v}(3;j) $ are given by Eq.(\ref{31}) and are as follows:
     $$
  l_{11}(t_{22})=1, \quad l_{11}(t_{11})=\cosh Jv, \quad
  \tilde{l}_{11}(\tilde{t}_{11})=-J^{-1}\sinh Jv, \quad
     $$
     $$
  l_{12}(\tilde{t}_{21})=-ij_1^2J^{-1}\sinh Jv,           \quad
  l_{12}(\tilde{t}_{12})=ij_1^2(2J)^{-1}(\sinh 3Jv/2 + \sinh Jv/2 ),
     $$
     $$
  l_{12}(t_{12})=(\cosh 3Jv/2 - \cosh Jv/2 )/2=\tilde{l}_{12}(\tilde{t}_{12}),
  \quad
\tilde{l}_{12}({t}_{21})=ij_2^2J^{-1}\sinh Jv,
     $$
     $$
  \tilde{l}_{12}(t_{12})=-ij_2^2(2J)^{-1}(\sinh 3Jv/2 + \sinh Jv/2 ), \quad
  l_{21}(\tilde{t}_{12})=-ij_1^2J^{-1}\sinh Jv,     \quad
     $$
     $$
  l_{21}(\tilde{t}_{21})=ij_1^2(2J)^{-1}(\sinh 3Jv/2 + \sinh Jv/2 ), \quad
   \tilde{l}_{21}(t_{12})=ij_2^2J^{-1}\sinh Jv,
    $$
    $$
 \tilde{l}_{21}(t_{21})=-ij_2^2(2J)^{-1}(\sinh 3Jv/2 + \sinh Jv/2 ), \quad
 l_{13}(t_{13})=(\cosh 2Jv - 1)/2 =\tilde{l}_{13}(\tilde{t}_{13}),  \quad
    $$
     $$
  l_{21}(t_{21})=(\cosh 3Jv/2 - \cosh Jv/2 )/2=\tilde{l}_{21}(\tilde{t}_{21}),
     $$
     \begin{equation}
 l_{13}(\tilde{t}_{13})=-iJ^{-1}(2\sinh Jv - \sinh 2Jv),             \quad
 \tilde{l}_{13}({t}_{13})=iJ(2\sinh Jv - \sinh 2Jv),
    \label{42}
    \end{equation}
where $ J=j_1j_2. $ Only nonzero expressions are written out above.
According with the additional relations (\ref{34}) there are three
independent generators of $ so_{v}(3;j), $ for example,
     $ l_{11}, l_{12}, \tilde{l}_{12}. $
Their commutation relations follow from Eq.(\ref{33})
$$
l_{11}l_{12} {\cosh}Jz - l_{12} l_{11} =
 l_{11} \tilde l_{12}
ij_1^{2}J^{-1}{\sinh}Jz,
$$
 $$
l_{11}  \tilde  l_{12}  {\cosh}Jz  - \tilde l_{12} l_{11} = -
l_{11} l_{12} ij_2^{2}J^{-1}{\sinh}Jz ,
$$
\begin{equation}
\left [ l_{12} , \tilde l_{12} \right ] = -
\left ( l_{11}^2  - 1 \right )
iJ{\sinh}Jz .
\end{equation}

The quantum analogue of the  universal enveloping
 algebra of CK algebra $ so(3;j)=\{X_{01},X_{02},X_{12}\}  $ with
 the rotation generator $ X_{02} $ as the primitive
element of the Hopf algebra has been given in\cite{Val--93}.
 The Hopf algebra structure of
 $ so_v(3;j;X_{02}) $  is given by
      $$
\Delta X_{02}  =  I\otimes X_{02}+X_{02}\otimes I,
      $$
      $$
\Delta X  =  e^{-vX_{02}/2}\otimes X+X\otimes e^{vX_{02}/2},
\quad  X=X_{01},X_{12},
      $$
      $$
\epsilon(X_{01})  =  \epsilon(X_{02})=\epsilon(X_{12})=0, \quad
 S(X_{02})  =  -X_{02},
      $$
      $$
 S(X_{01})  =  -X_{01}\cos{Jv/2}+
j_1^2X_{12}J^{-1}\sin{Jv/2},
      $$
      $$
 S(X_{12})  =  -X_{12}\cos{Jv/2}-
j_2^2X_{01}J^{-1}\sin{Jv/2},
      $$
 \begin{equation}
 [X_{01},X_{02}]=j_1^2X_{12},\quad   [X_{02},X_{12}]=j_2^2X_{01},\quad
 [X_{12},X_{01}]={ \sinh{ vX_{02}} \over {v}}.
\label{43}
 \end{equation}
The isomorphism of $ so_v(3;j;X_{02}) $ and quantum algebra
     $ so_v(3;j) $ is easily established with the help of the
following relations between generators
        $$
   l_{11}=e^{vX_{02}}, \quad
 j_1^{-1}l_{12}=j_2EX_{01}e^{vX_{02}/2}, \quad
j_2^{-1} \tilde{l}_{12}=-j_1EX_{12}e^{vX_{12}/2},
        $$
 \begin{equation}
E  =  i{\left (vJ^{-1}\sin Jv \right ) }^{1/2}e^{-Jv}.
\label{44}
\end{equation}
So the quantum orthogonal CK algebras may be constracted
both as the dual to the quantum group and by the contractions
of quantum orthogonal  algebras.

\section{ Conclusion }
\noindent


Contractions are the method of receiving a new Lie groups (algebras)
from the unitial ones, in particular, the nonsemisimple qroups
(algebras) from the simple ones. In the traditional approach{\cite{Wig--53}}
 this is achieved by introduction of a real zero
tending parameter $ \epsilon \rightarrow 0 $. In our approach{\cite{2}},
contractions are described by the dual valued parameters
 $ j_k. $
In the case of FRT theory of quantum
groups
these contractions supplemented with the
appropriate transformations of the deformation parameter
lead to realization of nonsemisimple quantum qroups as Hopf
 algebras of noncommutative functions with dual variables.
     So at least in CK scheme  from mathematical
     point of view the  contraction procedure is nothing else
as the replacement of complex number field $ {\bf C} $ by
     the  algebra  $ {\bf D}(\iota) $ for both classical and
     guantum groups.

\end{document}